\documentclass[11pt]{article}
\usepackage{amssymb, amsmath}
\def\refer#1{~\ref{#1}}
\def\refeq#1{~(\ref{#1})}
\def\ccite#1{~\cite{#1}}

\def\longformule#1#2{
\displaylines{ \qquad{#1} \hfill\cr \hfill {#2} \qquad\cr } }
\def\inte#1{
\displaystyle\mathop{#1\kern0pt}^\circ }

\let\lam=\lambda

\let\D=\Delta

\let\wt=\widetilde

\def\cB{{\cal B}}

\def\cD{{\cal D}}

\def\cI{{\cal I}}

\def\cL{{\cal L}}

\def\cR{{\cal R}}

\def\cT{{\cal T}}

\def\virgp{\raise 2pt\hbox{,}}
\def\cdotpv{\raise 2pt\hbox{;}}

\def\eqdefa{\buildrel\hbox{{\rm \footnotesize def}}\over =}
\def\Id{\mathop{\rm Id}\nolimits}

\def\C{\mathop{\bf C\kern 0pt}\nolimits}
\def\DD{\mathop{\bf D\kern 0pt}\nolimits}
\def\K{\mathop{\bf K\kern 0pt}\nolimits}
\def\N{\mathop{\bf N\kern 0pt}\nolimits}
\def\Q{\mathop{\bf Q\kern 0pt}\nolimits}
\def\R{\mathop{\bf R\kern 0pt}\nolimits}
\def\SS{\mathop{\bf S\kern 0pt}\nolimits}
\def\ZZ{\mathop{\bf Z\kern 0pt}\nolimits}
\def\TT{\mathop{\bf T\kern 0pt}\nolimits}
\def\P{\mathop{\bf P\kern 0pt}\nolimits}
\def\M{{\rm M}}
\def\St{{\rm S}}
\def\H{{\rm H}}
\newcommand{\ds}{\displaystyle}

\newcommand{\Z}{{\mathbf Z}}

\def\dive{\mathop{\rm div}\nolimits}

\def\Supp{\mathop{\rm Supp}\nolimits\ }

\newcommand{\ef}{ \hfill $ \blacksquare $ \vskip 3mm}
\newcommand{\beq}{\begin{equation}}
\newcommand{\eeq}{\end{equation}}
\newcommand{\ben}{\begin{eqnarray}}
\newcommand{\een}{\end{eqnarray}}
\newcommand{\beno}{\begin{eqnarray*}}
\newcommand{\eeno}{\end{eqnarray*}}
\newtheorem{defi}{Definition}[section]
\newtheorem{theo}{Theorem}
\newtheorem{lemma}{Lemma}[section]

\newtheorem{col}{Corollary}[section]
\newtheorem{prop}{Proposition}[section]

\begin{document}

\title{On  the global wellposedness  of the 3-D 
   Navier-Stokes equations with large initial data}
\author{Jean-Yves Chemin $ ^*$ and Isabelle  Gallagher $ ^\dag$\\
{\small $ ^*$ Laboratoire J.-L. Lions, Case 187 }\\
{\small Universit{\'e} Paris 6, 75252 Paris Cedex 05, FRANCE }\\
{\small chemin@ann.jussieu.fr }\\[2mm]  {\small $ ^\dag$
Institut de Math{\'e}matiques de Jussieu, Case 7012}\\
{\small Universit{\'e} Paris 7,
75251 Paris Cedex 05, FRANCE} \\
{\small Isabelle.Gallagher@math.jussieu.fr}}

\date{\today}
\maketitle

\begin{abstract}
We give a  condition for  the periodic, three dimensional, incompressible Navier-Stokes
 equations to be globally
wellposed. This condition is not a smallness condition on the initial
 data, as  
the data  is allowed to be arbitrarily
large in the    scale invariant 
space~$ B^{-1}_{\infty,\infty}$, which contains all the known
 spaces in which there is a global solution for small data. The
 smallness condition is rather a nonlinear type condition
on the initial data; 
  an explicit example of such initial data is constructed,  which is arbitrarily large and yet gives rise to a global, 
smooth solution.
\end{abstract}


\setcounter{equation}{0}
\section{Introduction}
The purpose of this text is to establish a condition of global wellposedness
for regular initial  data for the incompressible Navier-Stokes system on the
three dimensional torus~$\TT^3=(\R/2\pi\ZZ)^3$. Let us recall the system:
\[
(NS)\ \left\{
\begin{array}{c}
\partial_{t} u -\Delta u +u\cdot\nabla u=-\nabla p\\
\dive u =0\\
u_{t=0}=u_{0}.
\end{array}
\right.
\]
Here~$ u$ is a mean free three-component vector field~$ u = (u^1,u^2,u^3) =(u^h,u^3) $ representing the velocity 
of the fluid,  and~$
p$ is a scalar denoting the pressure, both are unknown functions of the space variable~$ x \in \TT^3$, and
the time variable~$ t \in \R^+$.  
We have chosen the kinematic viscosity of the fluid to be equal to one
for simplicity.
We recall that the pressure can be
eliminated by projecting~$ (NS) $ onto the space of divergence free
vector fields, using the Leray projector
$$
\P = \mbox{Id} - \nabla
\Delta^{-1} \mbox{div} .
$$ Thus we shall be using in the following the
equivalent system
$$
\partial_{t} u -\Delta u +\P (u\cdot\nabla u)=0. 
$$
Our motivation is the study of the size of the initial data yielding
global existence of solutions to that system, rather than
the minimal regularity one can assume on the initial data. Thus, in
all this work, we shall  assume that~$u_{0}$
is a mean free vector field with components in the Sobolev
space~$H^{\frac12}(\TT^3)$: we recall that~$H^{\frac12}(\TT^3)$ is a
scale invariant space for~$ (NS)$, and that smooth solutions exist for
a short time if the initial data belongs to~$H^{\frac12}(\TT^3)$,
globally in time if the data is small enough.  The problem
of global wellposedness for general data in~$H^{\frac12}(\TT^3)$ is
known to be open.   The search of smallness
conditions on~$u_{0}$ the least restrictive as possible is a long story, essentially
initiated by J. Leray (in the whole space~$\R^3$ but the phenomenon is similar in
the torus) in the seminal paper\ccite{lerayns}, continued in particular
by H. Fujita and T. Kato in\ccite{fujitakato},
M. Cannone, Y. Meyer et F. Planchon in\ccite{cannonemeyerplanchon}
and H. Koch and D. Tataru in\ccite{kochtataru}. The theorem proved
in\ccite{kochtataru}  claims that if~$\|u_{0}\|_{\partial BMO}$ is small, which
means that the components of~$u_{0}$ are derivatives of~$BMO$
functions and are small enough, then~$(NS)$ is globally wellposed
in the sense that a global (and of course unique) solution exists
in~$C(\R^+;H^\frac12)$. Our aim is to prove a theorem of global wellposedness
which allows for very large data in~$\partial BMO$,  under a nonlinear smallness
condition on the initial data. In fact the initial data will even be
large in~$ B^{-1}_{\infty,\infty}$,   which contains
strictly~$\partial BMO$ and which
is the largest scale invariant Banach space in which one can hope to
prove a wellposedness result.  Before stating the result, let us recall
that the question is only  meaningful in three or
more space dimensions. We recall indeed that according to
J. Leray~\cite{leray2d}, there is a unique, global
 solution to the two dimensional Navier-Stokes system as soon as the
 initial data is  in~$L^{2}(\TT^{2}), 
$ and if there is a forcing term it should belong for instance
to~$L^{1}(\R^{+}; L^{2}(\TT^{2}))$.

\medskip

\noindent In order to state our result, we shall need the following notation: 
one can decompose  any
function~$ f$ defined on~$ \TT^3$ as
$$
f = \overline f + \widetilde f, \quad \mbox{where} \quad \overline
f(x_1,x_2 ) =  \frac1{2\pi} \int_0^{2\pi} f(x_1,x_2,x_3) \: dx_3. 
$$
Similarly we shall define the horizontal mean~$ \overline u$ of any
vector field as~$ \overline u = (\overline u^1, \overline u^2,
\overline u^3)$. It will also be convenient to use the following
alternative notation: we denote by~$ \M$ the projector onto vector
fields defined on~$ \TT^2$, 
$$
\M f = \overline f \quad \mbox{and} \quad (\Id - \M)f = \widetilde f.
$$
we shall denote the heat semiflow by~$
\St(t)  = e^{t\Delta}
$. Finally let us define negative index Besov spaces.
\begin{defi}
\label{besov}
{\sl
Let~$s$ be a positive real number, and let~$p$ and~$q$  be two real
numbers in~$[1,+\infty]$. The Besov space~$ B^{-s}_{p,r}(\TT^{3})$
is the space of distributions in~$ \TT^3$ such that
$$
\|u\|_{ B^{-s}_{p,q}} \eqdefa \left\| t^{\frac{s}{2}}\|\St(t)  u\|_{L^{p}}\right\|_{L^{q}(\R^{+},
\frac{dt}{t})} < +\infty.
$$
}
\end{defi}
{\bf Remark }  we shall see an equivalent definition in terms of
Littlewood-Paley theory in Section~\ref{notation} (see Definition~\ref{besovlp}).

\medskip
\noindent
Now let us consider the following subspace of~$H^{\frac12}(\TT^{3})$, where we have noted, for all
vector fields~$a$ and~$b$,
$$
Q(a,b) \eqdefa \P \dive (a \otimes b + b \otimes a). 
$$
\begin{defi}
\label{definI(A,B)}{\sl
Let~$ A$ and~$ B$ be two positive real numbers and let~$p$ in~$]3,+\infty]$. We 
define the set
$$
\cI_p(A,B) = \left\{u_{0} \in H^{\frac 1 2}(\TT^3) \: \Big / \: \dive u_{0} = 0 \:\:
  \mbox{and} \:\: (\H1) ,\: (\H2) , \:(\H3)\:
  \mbox{are satisfied}\right\}, \quad \mbox{where}
$$
\beno
\displaystyle (\H1) \quad    &\ds \|\overline u_{0}\|_{L^2(\TT^2) } + 
\|\M {\bf P}  (  u_{F}\cdot \nabla   u_{F} )\|_{L^1(\R^+;L^2(\TT^2))} \leq A & \\
\displaystyle (\H2)  \quad    &\ds \|\widetilde{u_{0}}\|_{
  B^{-1}_{\infty,2}} \leq A&\\
\displaystyle (\H3)  \quad        &\ds\left\|
(\Id - \M) {\bf P} ( u_{F}\cdot \nabla u_{F} )
  +Q  ( u_{2D} ,  u_{F} )  
 \right\|_{L^1(\R^+; B^{-1+\frac3p}_{p,2})}   \leq B ,&
\eeno
where we have noted~$
u_{F} (t) = \St (t) \wt u_{0} $ and where~$u_{2D}$ is a three component vector field defined on~$\TT^{2}$, satisfying
the following two dimensional Navier-Stokes equation, in the case when
the initial
data is~$v _{0} = \overline u_{0} $ and the force is~$ f = -\M {\bf P}
(u_{F} \cdot \nabla u_{F}) $:
$$
(NS2D) \quad \left\{\begin{array}{c}
\partial_{t }  v  +  \P ( v^h  \cdot \nabla^{h}  v  ) - \Delta_{h}  v  
= f  \\
 v _{|t = 0} =  v_{0} ,
\end{array}
\right.
$$
 where~$\Delta_{h}$ denotes the horizontal 
 Laplacian~$\Delta_{h } = \partial_{1}^2 + \partial_{2}^{2}$ and
 where~$  \nabla^{h} = (\partial_1,\partial_2)$. 
}
\end{defi}
Now let us state the main result of this paper.
\begin{theo}\label{mainresult}
{\sl 
Let~$ p \in \: ]6,+\infty[$ be given. 
There is a constant~$C_{0} > 0$ such that the following
holds. Consider two positive real numbers~$A$
and~$B$ satisfying
\begin{equation}\label{smallness}
B \exp \left( C_0 A^2 \left(1 + A \log (e + A)\right)^2\right)\leq C_{0}^{-1}.
\end{equation}
Then for any vector field~$u_{0} \in \cI_p(A,B)$, there is a unique, global solution~$u$ to~$(NS)$ associated
with~$u_{0}$, satisfying
\[
u\in C_{b}(\R^{+};H^{\frac 1 2}(\TT^{3}))
\cap L^2(\R^+;H^{\frac 3 2}(\TT^3)).
\] 
}
\end{theo}
{\bf Remarks } 

\noindent 1) Condition~(\ref{smallness}) appearing in the statement of
 Theorem~\ref{mainresult}
should be understood as a nonlinear smallness condition on the initial data: the parameter~$A$, 
measuring through~$(H1)$ and~$(H2)$ the norm of the initial data in a scale-invariant space, may be as large
as wanted, as long as the parameter~$B$, which measures a nonlinear quantity in a scale-invariant 
space, is small enough. We give below an example of such initial data, which  is a smooth vector 
field
with   arbitrarily large~$ B^{-1}_{\infty,\infty}$ norm, and which generates a unique, 
global solution to~$(NS)$: see the statement of Theorem~\ref{example}. 
 
\medskip

\noindent  2) Some results of global existence for large data can be found in 
the literature. To our knowledge they all involve either an initial vector field 
which is close enough to a two dimensional
vector field (see for instance\ccite{raugelsellthindom},~\cite{gimrn} or~\cite{iftimie}), or initial 
data such that after a change of coordinates, the equation is transformed into 
the three dimensional rotating  fluid equations (for which global existence is 
known), see~\cite{bmn}. Here we are in neither of those situations. 

\bigskip

\noindent Let us now give an example where condition~(\ref{smallness}) holds. As mentioned in the 
remarks above, in that example the initial data can be arbitrarily large in~$ B^{-1}_{\infty,\infty}$, and nevertheless
generates a global solution. We have noted by~$\widehat u$ the Fourier transform of any vector
field~$u$. 
\begin{theo}\label{example}
{\sl 
Let~$N_{0}$ be a given positive integer. A positive integer~$N_{1}$ exists 
such, if~$N$ is an integer larger than~$N_{1}$, it satisfies  the following
properties. If~$v_{0}^h$ is any two component,  
divergence free vector fields defined on~$\TT^{2}$  such that
$$
\Supp\widehat v_{0}^h \subset [-N_{0},N_{0}]^{2}
\quad\hbox{and}\quad
\|v^h_{0}\|_{L^2(\TT^2)} \leq(\log N)^{\frac 1 9},
$$
then a unique, global smooth solution  to~$(NS)$ exists,  associated with the initial data
$$
u_{0}(x) = \left(
N v_{0}^h(x_{h}) \cos (Nx_{3}) , - \dive_{h}  v_{0}^h(x_{h})  \sin(Nx_{3})
\right).
$$
Moreover the vector field~$u_{0}$ satisfies
\begin{equation}\label{largeC-1}
\|u_{0}^{h}\|_{ B^{-1}_{\infty,\infty}} \geq \frac{1}{4\pi \sqrt e} \|v_{0}^h\|_{L^{2}(\TT^{2})} .
\end{equation}
}
\end{theo}
{\bf Remarks }

\noindent 1)  Since the~$L^{2} $ norm of~$v_{0}^h$ can be chosen arbitarily
large,  the lower bound given in~(\ref{largeC-1}) implies that the~$ B^{-1}_{
\infty,\infty}$ norm of the initial data may be chosen arbitrarily large. 

\medskip

\noindent 2) One can rewrite this example in terms of the Reynolds
number of the fluid: let~${\rm re} \in \N$ be its Reynolds
number, and define the rescaled velocity field~$\displaystyle v(t,x ) = \frac1{\rm
  re} u (\frac{t}{\rm  re}\virgp x )$.  Then~$ v$ satisfies the
Navier-Stokes equation
$$
\partial_t v + {\bf P} (v\cdot \nabla v) - \nu  \Delta v= 0
$$
where~$ \nu = 1 / {\rm re}$, 
and Theorem~\ref{example} states the following: if~$ v_{|t = 0}$ is
equal to
$$
v_{0, \nu } = \left(
v^h_{0} (x_h) \cos \left(\frac{  x_3}{\nu}\right) , - \nu
\dive_h v^h_{0} (x_h) \sin\left(\frac{  x_3}{\nu}\right)
\right)
$$
where~$\widehat v^h_{0}  $ is supported in~$ [-N_0,N_0]^2$
and satisfies
$$
 \|v_{0}^{h}\|_{L^{2}(\TT^{2})} \leq
 \Bigl (\log \frac1\nu\Bigr)^{\frac 1 9}\virgp
$$
then for~$\nu $ small enough there is a unique, global, smooth
solution. 

\bigskip

\noindent The rest of the paper is devoted to the proof of Theorems~\ref{mainresult} and~\ref{example}. 
The proof of  Theorem~\ref{mainresult} relies on the following idea: if~$u$ denotes the solution
of~$(NS)$ associated with~$u_{0}$, which exists at least for a short
time since~$ u_0$ belongs to~$ H^{\frac12}(\TT^3)$, then it can be decomposed
as follows, with the notation of Definition~\ref{definI(A,B)}:
$$
u = u^{(0)}  + R, \quad \mbox{where} 
\quad u^{(0)} = u_{F}+u_{2D}.
$$

\noindent Note that the Leray theorem in dimension two mentioned
above, namely
 the existence and uniqueness of a smooth solution for an initial data in~$L^{2}(\TT^{2})$ and a forcing
term in~$L^{1}(\R^{+};L^{2}(\TT^{2}))$, 
holds even if the vector fields have three components rather than two
 (as is the case for the equation satisfied by~$u_{2D} $).  One
 notices that the vector field~$R$  satisfies the perturbed
 Navier-Stokes system
$$
(PNS) \quad \left\{\begin{array}{c}
\partial_{t} R + \P (R \cdot \nabla R) + Q (u^{(0)},R) - \Delta R = F \\
R_{|t = 0} = R_{0} ,
\end{array}
\right.
$$
where
$$
R_{|t = 0} = 0\quad \mbox{and} \quad F =  - (\Id - \M) {\bf P}(u_{F} \cdot \nabla u_{F}) 
- Q (u_{F} , u_{2D} ) .
$$
The proof of Theorem~\ref{mainresult} consists in studying both systems, the 
two dimensional Navier-Stokes
system and the perturbed three dimensional Navier-Stokes system. In particular 
a result on the two dimensional Navier-Stokes
system will be proved in Section\refer{NS2DL2Linfty}, which, as far as we 
know  is new, and may have its own interest. It is stated below.
\begin{theo}
\label{NS2DLinfty}
{\sl There is a constant~$ C>0$ such that the following result holds. 
Let~$v$ be the solution of~$(NS2D)$ with initial data~$v_{0}\in L^2$
and external force~$ f$ in~$L^1 (\R^+;L^2)$.  Then we have
\[
\|v\|_{L^2(\R^+;L^\infty)}^2 \leq C  E_{0}
\Bigl(1+E_{0}\log^2 (e+E^{\frac 1 2}_{0})\Bigr)
\quad\hbox{with}\quad
E_{0}\eqdefa
\|v_{0}\|_{L^2}^2+\Bigl(\int_{0}^\infty\|f(t)\|_{L^2}dt\Bigr)^2.
\]}
\end{theo}
\noindent The key to the proof of Theorem~\ref{mainresult} is the proof of the global
wellposedness of the perturbed three dimensional system~$ (PNS)$. That 
is achieved in Section~\ref{PNS3Dglobal} below, where a general
statement is proved, concerning the global wellposedness of~$ (PNS)$
for general~$ R_0$ and~$ F$ satisfying a smallness condition. That
result is joint to Theorem~\ref{NS2DLinfty} to prove
Theorem~\ref{mainresult} in Section~\ref{proof}. 
Finally Theorem~\ref{example} is proved
in Section~\ref{exampleproof}. The coming section is devoted to some
notation and the recollection of well known results on Besov spaces
and the Littlewood-Paley theory which will be used in the course of
the proofs.

\section{Notation and useful results on Littlewood-Paley theory}\label{notation}
\setcounter{equation}{0}
In this short section we shall present some well known facts on the
Littlewood-Paley theory. Let us start by giving the definition of
Littlewood-Paley operators on~$ \TT^d$. 
\begin{defi}
\label{d1}
{\sl
Let $\chi$ be a nonnegative function in $C^\infty(\TT^{d})$ such that $\widehat\chi =
1$ for~$|\xi|\leq 1$ and~$\widehat\chi= 0$ for~$|\xi|>2$, and define
$
\chi_{j}(x)= 2^{jd}\chi(2^{j}|x|).$
Then the Littlewood-Paley frequency localization operators are defined by
$$ 
S_{j} = \chi_{j}\ast\cdot, \quad 
\Delta_{j} = S_{j} - S_{j-1}.
$$ 
}
\end{defi}
As is well known, one of the interests of this decomposition is that
the~$ \Delta_j$ operators allow to count derivatives easily. More precisely we recall
the  Bernstein inequality. A constant~$C$ exists such that
\begin{equation}
\label{bernstein}
\forall \: k \in \N, \:
\forall  \: 1 \leq p \leq q \leq \infty, \quad 
\sup_{|\alpha| = k}\|\partial^\alpha\Delta_j u\|_{L^q(\TT^d)}  
\leq C^{k+1}2^{jk}2^{jd \left(\frac{1}{p} -  \frac{1}{q}\right )}
\|\Delta_j u\|_{L^p (\TT^d)} .
\end{equation}
 Using those operators we can give a definition of Besov spaces for all
indexes, and we recall the classical fact that the definition in the
case of a negative index coincides with the definition given in the
introduction using the heat kernel (Definition~\ref{besov} above). 
\begin{defi}
\label{besovlp}
{\sl 
Let~$f$ be a mean free function in~$\cD'(\TT^{d})$, and let~$ s \in
\R$ and~$ (p,q) \in [1,+\infty]^2$ be given real numbers. Then~$f$ belongs to the Besov
space~$ B^{s}_{p,q}(\TT^d)$ if and only if
$$
\|f\|_{ B^{s}_{p,q}} \eqdefa  \left\|2^{js}\| \Delta_{j}
f\|_{L^{p}}\right\|_{\ell^{q}(\Z)} < +\infty.
$$
}
\end{defi}
Using the Bernstein inequality~(\ref{bernstein}), it is easy to see that the
 following continuous embedding
 holds:
\begin{equation}\label{embedding}
 B^{s + \frac{d}{p_1}}_{p_1,r_1}
  (\TT^d) \hookrightarrow  B^{s + \frac{d}{p_2}}_{p_2,r_2}
  (\TT^d) ,
\end{equation}
for all real
 numbers~$ s,p_1,p_2,r_1,r_2$
 such that~$ p_i$ and~$ r_i$ 
belong to the interval~$[1,\infty] $ and such that~$ p_1 \leq p_2 $
 and~$ r_1 \leq r_2 $. 

\medskip

\noindent We recall that Sobolev spaces are special cases of Besov
spaces, since~$ H^s = B^s_{2,2}$. 

\medskip

\noindent Throughout this article we shall denote by the letters~$ C$ or~$ c$ all
universal constants. we shall sometimes replace an inequality of the
type~$f \leq C g$ by~$ f \lesssim g$. we shall also denote by~$
(c_j)_{j \in \ZZ}$
any sequence of norm~1 in~$ \ell^2(\ZZ)$.

\section{An~$L^\infty$ estimate for Leray solutions in dimension two}
\label{NS2DL2Linfty}
\setcounter{equation}{0}
The purpose of this section is the proof of Theorem\refer{NS2DLinfty}.
Let us write  the solution~$v$ of~$(NS2D)$ as the sum of~$v_{1}$ 
 and~$v_{2}$  with
 \beq
 \label{NS2DL2Linftyeq1}
 \left\{
 \begin{array}{c}
 \partial_{t}v_{1}-\Delta_{h} v_{1} = {\bf P}f\\
 {v_{1}}_{|t=0} = v_{0}
 \end{array}
 \right.
 \quad\hbox{and}\quad
 \left\{
 \begin{array}{c}
 \partial_{t}v_{2}-\Delta_{h} v_{2} =  - {\bf P}\dive (v\otimes v)\\
 {v_{2}}_{|t=0} = 0.
 \end{array}
 \right.
 \eeq
Duhamel's formula gives
\[
v_{1}(t) = e^{t\Delta} v_{0}+\int_{0}^t e^{(t-t')\Delta}{\bf P}
f(t')dt', 
\]
thus we get that
\[
\|v_{1}\|_{L^2(\R^+;L^\infty)} \leq 
\|e^{t\Delta} v_{0}\|_{L^2(\R^+;L^\infty)} 
+\int_{0}^\infty 
\left\|e^{\tau\D}f(t)\right\|_{L^2(\R^+_{\tau};L^\infty)}dt.
\]
Due to~(\ref{embedding}), we have $L^2\hookrightarrow
B^{-1}_{\infty,2}$ so by Definition~\ref{besov} we get that
\beq
\label{NS2DL2Linftyeq2}
\|v_{1}\|_{L^2(\R^+;L^\infty)} \lesssim
\|v_{0}\|_{L^2} 
+\int_{0}^\infty \|f(t)\|_{L^2}\:dt.
\eeq
Now let us estimate~$\|v_{2}\|_{L^2(\R^+;L^\infty)}$. It relies on the 
following technical proposition.
\begin{prop}
\label{estimnormtilde2D}
{\sl
Let~$v$ be the solution of~$(NS2D)$ with initial data~$v_{0}$ in~$L^2$ 
and  external force~$ f$ in~$L^1 (\R^+;L^2)$.  Then we have
\[
\sum_{j}\|\D_{j}v\|^2_{L^\infty(\R^+;L^2)}\lesssim 
E_{0}(e+E_{0}^{\frac 1 2})\quad\hbox{with}\quad
E_{0} =
\|v_{0}\|_{L^2}^2+\Bigl(\int_{0}^\infty\|f(t)\|_{L^2}dt\Bigr)^2.
\]
}
\end{prop}

\noindent{\bf Proof of Proposition\refer{estimnormtilde2D} }
Applying~$\Delta_{j}$ to the~$(NS2D)$ system and doing an~$L^2$ energy 
estimate gives, neglecting (only here) the smoothing effect of the heat flow,
\[
\|\D_{j}v (t)\|_{L^2}^2 \leq 
\|\D_{j}v_{0}\|_{L^2}^2
+\int_{0}^t \left|\left(\D_{j}(v(t')\cdot\nabla v(t'))|\D_{j}v(t')\right)_{L^2}\right|dt'
+ \int_{0}^t\left| \langle \D_{j}f(t'),\D_{j}v(t')\rangle  \right| dt'.
\]
Lemma 1.1  of\ccite{chemin10} and the conservation of energy tell us that
\beno
\left|\left(\D_{j}(v(t)\cdot\nabla v(t))|\D_{j}v(t)\right)_{L^2}\right| 
& \lesssim & 
c_{j}(t)\|\nabla v(t)\|_{L^2}\|v(t)\|_{L^2} 2^j \|\D_{j}v(t)\|_{L^2}\\
& \lesssim  & c_{j}^2(t) \|\nabla v(t)\|_{L^2}^2\|v(t)\|_{L^2}\\
& \lesssim & E_{0}^{\frac 1 2}c_{j}^2(t) \|\nabla v(t)\|_{L^2}^2.
\eeno
Since
\beno
\left| \langle \D_{j}f(t),\D_{j}v(t)\rangle  \right| & \leq  &
\|\D_{j}f(t)\|_{L^2} \|\D_{j}v(t)\|_{L^2} \\
& \lesssim &  E_{0}^{\frac 1 2}c_{j}^2(t) \|f(t)\|_{L^2}, 
\eeno
 we infer that
 \[
\|\D_{j}v \|_{L^\infty(\R^+;L^2)}^2  \lesssim 
  \|\D_{j}v_{0}\|_{L^2}^2 +E_{0}^{\frac 1 2}
\int_{0}^\infty c^2_{j}(t)\left(\|\nabla v(t)\|_{L^2}^2
+\|f(t)\|_{L^2}\right)dt.
\]
Taking the sum over $j$ concludes the proof of the proposition.\ef
 
\medbreak\noindent
 {\bf Conclusion of the proof of  Theorem\refer{NS2DLinfty} }
Let us first observe that interpolating the result of
 Proposition~\ref{estimnormtilde2D}  with the energy estimate, we find
 that a 
constant~$C$ exists such that, for any~$p$ in~$[2,\infty]$, we have
\beq
\label{espacetildeinterpoleq}
\sum_{j}2^{j\frac 4 p }\|\D_{j}v\|^2_{L^p(\R^+;L^2)} \leq C E_{0}
(e+E_{0}^{\frac 1 2})^{1-\frac 2 p}.
\eeq
Then by  Bernstein's inequality~(\ref{bernstein}), we have
\[
2^{-j\left(1-\frac 2 p\right)}\|S_{j}v\|_{L^p(\R^+;L^{\infty})}
\leq C \sum_{j'\leq j-1} 2^{(j'-j)\left(1-\frac 2 p\right)}
2^{j'\frac 2 p} \|\D_{j'}v\|_{L^p(\R^+;L^{2})} .
\]
 Using Young's inequality on series and\refeq{espacetildeinterpoleq}, 
 we infer that a constant~$C$ exists such that,  for any~$p$ in~$]2,\infty]$,
\beq
\label{NS2DL2Linftyeq11}
2^{-j\left(1-\frac 2  p\right)}\|S_{j}v\|_{L^p(\R^+;L^{\infty})}
\leq C c_{j}\frac p {p-2} E_{0}^\frac12 (e+E_{0}^{\frac 1 2})^{\frac12-\frac 1 p}.
\eeq
Now using Bernstein's inequality and Fourier-Plancherel, we get by~(\ref{NS2DL2Linftyeq1})
 \ben
\| \D_{j} v_{2}(t)\|_{L^\infty} & \lesssim & 
2^{j}\| \D_{j} v_{2}(t)\|_{L^2}\nonumber\\
& \lesssim & \label{NS2DL2Linftyeq3}
 2^{2j}\int_{0}^t e^{-c2^{2j}(t-t')}
\left\|\D_{j} {\bf P} \left (v(t')\otimes v(t')\right)\right\|_{L^2}dt'.
 \een
Using Bony's decomposition, let us write that for any~$ a$ and~$ b$, 
\[
\D_{j}(a(t)b(t)) = \sum_{j'\geq j-N_{0}}\D_{j}( S_{j'}a(t) \D_{j'}b(t))
+ \sum_{j'\geq j-N_{0}} \D_{j}(\D_{j'}a(t) S_{j'+1}b(t)).
\]
We have
\[
\|S_{j'}a \D_{j'}b\|_{L^{\frac {2p} {p+2}}(\R^+;L^2)}\leq
\|S_{j'}a\|_{L^p(\R^+;L^\infty)}\|\D_{j'}b\|_{L^2(\R^+;L^2)}.
\]
Using\refeq{NS2DL2Linftyeq11}, we deduce that  a constant~$C$ exists 
such that,  for any~$p$ in~$]2,\infty]$,
\beno
\left\|\D_{j} {\bf P}  \left (v\otimes v\right)\right
\|_{L^{\frac {2p}{p+2}}(\R^+;L^2)}
& \leq & C  \sum_{j'\geq j-N_{0}} \|S_{j'}v\|_{L^p(\R^+;L^\infty)}
\|\D_{j'}v\|_{L^2(\R^+;L^2)}\\
& \leq  & C    \frac p {p-2} E_{0}^\frac12 (e+E_{0}^{\frac 1
  2})^{\frac12-\frac 1 p} \!
\sum_{j'\geq j-N_{0}} c_{j'}\|\D_{j'} v\|_{L^2(\R^+;L^2)} 
2^{-j'\left(\frac2p - 1\right)}.
\eeno
Using Young's inequality in time in\refeq{NS2DL2Linftyeq3} gives
\beno 
\| \D_{j} v_{2}\|_{L^2(\R^+;L^\infty)}  & \leq &   C  2^{2j}
\|e^{-c2^{2j}\cdot}\|_{L^{\frac p {p-1}}}
\left\|\D_{j} {\bf P}\left (v\otimes v\right)\right
\|_{L^{\frac {2p}{p+2}}(\R^+;L^2)}\\
& \leq  & C \frac p {p-2} E_{0}^\frac12 (e+E_{0}^{\frac 1
  2})^{\frac12-\frac 1 p}  \sum_{j'\geq j-N_{0}} c_{j'}
 \|\D_{j'} v\|_{L^2(\R^+;L^2)} 2^{j'} 2^{(j-j')\frac2p }.
\eeno
By Young's inequality on series we find that a constant~$C$ exists such that, for any~$p$ 
in~$]2,\infty[$, 
\[
\| \D_{j} v_{2}\|_{L^2(\R^+;L^\infty)} \leq 
C c_{j}^2 \frac {p^2} {p-2} E_{0}^\frac12 (e+E_{0}^{\frac 1
  2})^{\frac12-\frac 1 p}
\]
and thus
\[
\|v_{2}\|_{L^2(\R^+;L^\infty)} \leq C \frac {p^2} {p-2} E_{0} (e+E_{0}^{\frac 1
  2})^{\frac12-\frac 1 p}.
\]
Then let us choose~$p$ such that
\[
\frac 2 p = 1-\frac 2 {\log (e+E_{0}^{\frac 1 2})}\cdotp
\]
Then we have that
\beq\label{estimatev2}
\|v_{2}\|_{L^2(\R^+;L^\infty)} \leq C E_{0}\log (e+E_{0}^{\frac
    1 2}),
\eeq
and putting~(\ref{NS2DL2Linftyeq2}) and~(\ref{estimatev2}) together proves Theorem\refer{NS2DLinfty}.\ef

\medskip

\noindent This theorem will enable us to infer the following useful corollary.
\begin{col}\label{u0L2linfty}
{\sl
Let~$ p \in \: ]2,+\infty[$ and let~$ u_0$ be a vector field in~$ {\mathcal I}_p(A,B)$. Then~$ u^{(0)} = u_F +
u_{2D}$ satisfies
$$
\|u^{(0)}\|_{L^2 (\R^+;L^\infty)}^2 \lesssim 
A^2 (1 + A\log (e + A))^2.
$$
}
\end{col}
{\bf Proof of Corollary~\ref{u0L2linfty} } As  in the proof of~(\ref{NS2DL2Linftyeq2}) above, we have clearly by
Definition~\ref{besov} and~$ (\H2)$,
\beno
\|u_F\|_{L^2 (\R^+;L^\infty)} &\leq& \|\widetilde
u_0\|_{B^{-1}_{\infty,2}} \\
& \leq &A .
\eeno 
Then by Theorem\refer{NS2DLinfty} we have
$$
\|u_{2D}\|_{L^2 (\R^+;L^\infty)}^2 \lesssim 
E_0 (1 +  E_0 \log ^2(e +  E_0^\frac12 )),
$$
where by definition of~$ E_0$ and by~$ (\H1)$, 
\beno
E_0& = &\|\overline u_0\|_{L^2}^2 + \|\M {\bf P} (u_F \cdot \nabla  u_F)\|_{L^1
  (\R^+;L^2)}^2 \\
&\leq &A^2 .
\eeno
As a result we get
$$
\|u^{(0)}\|_{L^2 (\R^+;L^\infty)}^2 \lesssim  A^2 + A^2 (1 + A
\log (e + A))^2
$$
and the corollary is proved. 
 \ef

\section{Global wellposedness of the perturbed system}
\label{PNS3Dglobal}
\setcounter{equation}{0}
In this section we shall study the  global wellposedness of the
system~$(PNS)$. The result is the following.
\begin{theo}\label{PNSglobal}
{\sl 
Let~$ p \in \: ]3,+\infty[$ be given. There is a constant~$C_{0} > 0$
such that for any~$ R_0$ in~$
B^{-1+\frac3p}_{p,2}$,  $ F$ in~$ L^1(\R^+;B^{-1+\frac3p}_{p,2})$
and~$u^{(0)} $ in~$ L^2(\R^+;L^\infty)$ satisfying
\begin{equation}\label{smallnessPNS}
\|R_0\|_{B^{-1+\frac3p}_{p,2}} +
\|F\|_{L^1(\R^+;B^{-1+\frac3p}_{p,2})} \leq C_0^{-1} e^{-C_0 \|u^{(0)}\|^2_{L^2(\R^+;L^\infty)}},
\end{equation}
 there is a unique, global solution~$R$ to~$(PNS)$ associated
with~$R_{0}$ and~$ F$, such that
\[
R\in C_{b}(\R^{+};B^{-1+\frac3p}_{p,2}) 
\cap L^2(\R^+;B^{\frac3p}_{p,2}).
\] 
}
\end{theo}
{\bf Proof of Theorem~\ref{PNSglobal}     }
Using Duhamel's formula, the system~$(PNS)$ turns out be
\beno
R & = &\cR_{0}+  L_{0}R+B_{NS}(R,R)\quad\hbox{with}\\
\cR_{0}(t) & \eqdefa &e^{t\Delta} R_0 +\int_{0}^t e^{(t-t')\Delta}F(t')dt'\,,\\
L_{0}R (t) & \eqdefa  & -\int_{0}^t e^{(t-t')\Delta} Q(u^{(0)}(t'), R(t'))dt'\quad\hbox{and}\\
B_{NS}(R,R)(t) & \eqdefa & -\int_{0}^t e^{(t-t')\Delta}{\bf P}
\dive\left(R(t')\otimes R(t')\right)dt'.
\eeno
The proof of the global wellposedness of~$ (PNS)$
 relies on the following classical fixed point lemma in a Banach space, the 
proof of which is omitted. 
\begin{lemma}
\label{cacciopoli+}
{\sl Let~$X$ be a Banach space, let~$L$ be a  continuous linear map 
from~$X$ to~$X$,  and let~$B$ be a bilinear map from~$X\times X$ to~$X$. 
Let us define
\[
\|L\|_{\cL(X)} \eqdefa \sup_{\|x\|=1} \|Lx\|\quad\hbox{and}\quad
\|B\|_{\cB(X)} \eqdefa \sup_{\|x\|=\|y\|=1} \|B(x,y)\|.
\]
If~$\|L\|_{\cL(X)}<1$, then for any~$x_{0}$ in~$X$ such that
\[
\|x_{0}\|_{X}< \frac{ (1-\|L\|_{\cL(X)})^2} {4\|B\|_{\cB(X)} }\virgp
\]
the equation
\[
x=x_{0}+Lx+B(x,x)
\]
has a unique solution in the ball of center~$0$ and 
radius~$\ds \frac {1-\|L\|_{\cL(X)}}{2\|B\|_{\cB(X)}}\cdotp$}
\end{lemma}
Solving system~$(PNS)$ consists therefore in finding a space~$X$ in  which 
we shall be able to apply Lemma\refer{cacciopoli+}. Let us define, for
any positive real number~$\lam$ and for any~$p$ in~$]3,\infty]$, 
the following space.
\begin{defi}
\label{definespacetildeapoids}
{\sl The space~$X_{\lam}$ is the space of distributions~$a$ on~$\R^+\times \TT^3$
such that
\[
\displaylines{
\|a\|^2_{X_{\lam}}\eqdefa \sum_{j} 2^{-2j\left(1-\frac 3 p\right)}
\Bigl(\|\Delta _{j}a_{\lam}\|^2_{L^\infty(\R^+;L^p)}
+2^{2j}\|\Delta _{j}a_{\lam}\|^2_{L^2(\R^+;L^p)}\Bigr)<\infty
\quad\hbox{with}\cr
a_{\lam}(t)\eqdefa 
\exp\Bigl(-\lam\int_{0}^t\|u^{(0)}(t')\|_{L^\infty}^2dt'\Bigr)a(t).
}
\]}
\end{defi}

\noindent{\bf Remark } If~$a$ belongs to~$X_{\lam}$, then~$a_{\lam}$ 
belongs to~$L^\infty(\R^+; B^{-1+\frac 3 p}_{p,2})\cap 
L^2(\R^+;B^{\frac 3 p}_{p,2})$ and, as~$u^{(0)}$ is 
in~$L^2(\R^{+};L^\infty)$, we have
\[
\|a\|_{L^\infty(\R^+; B^{-1+\frac 3 p}_{p,2})} 
+\|a\|_{L^2(\R^+; B^{\frac 3 p}_{p,2})} \leq 
\|a\|_{X_\lam}\exp \left(\lam\|u^{(0)}\|^2_{L^2(\R^+;L^\infty)}\right).
\]
The fact that~$X_{\lam}$ equipped with this norm is a Banach space is 
a routine exercise left to the reader.  The introduction of this space is 
justified by the following proposition which we shall prove at the end
of this section.
\begin{prop}
\label{propestimconjugaison}
{\sl For any~$p$ in~$]3,\infty[$, a constant~$C$ exists such that, for any 
positive~$\lam$,
\[
\| L_{0}\|_{  \cL(X_{\lam})} \leq \frac C {\lam^{\frac 1 2} }
\quad \hbox{and}\quad
\|B_{NS}\|_{\cB(X_{\lam})} \leq C e^{\lam \|u^{(0)}\|^2_{L^2(\R^+;L^\infty)}}.
\]
}
\end{prop}

\noindent
{\bf Conclusion of  the proof  of Theorem\refer{mainresult}  } In order to apply Lemma\refer{cacciopoli+}, 
let us choose~$\lam$ such that~$\|L_{0}\|_{ \cL(X_{\lam})}\leq 1/2$. 
Then, the condition required to apply Lemma\refer{cacciopoli+} is
\beq
\label{estimF0}
\|\cR_{0}\|_{X_{\lam}} \leq \frac{1}{16C}  e^{-4 C^2 \|u^{(0)}\|^2_{L^2(\R^+;L^\infty)}}.
\eeq
In order to ensure this condition, let us recall  Lemma~2.1 of\ccite{chemin20}.
\begin{lemma}
\label{regulheatflow}
{\sl A constant~$c$ exists such that, for any integer~$j$, any positive real number~$t$ 
and any~$p$ in~$[1,\infty]$,
\[
\|\D_{j}e^{t\Delta} a \|_{L^p} \leq 
\frac 1 c e^{-c2^{2j}t}\|\Delta_{j}a\|_{L^p}.
\]}
\end{lemma}
This lemma and the Cauchy-Schwarz inequality for the 
measure~$\|F(t')\|_{B^{-1+\frac 3 p}_{p,2}}dt'$ give
\beno
\|\Delta_{j}\cR_{0,\lam}(t)\|_{L^p} & \leq & C e^{-c2^{2j} t}
\|\Delta_{j} R_0\|_{L^p} + 
C \int_{0}^t e^{-c2^{2j}(t-t')} \|\Delta_{j}F(t')\|_{L^p}dt'\\
& \leq & C  2^{j\left(1-\frac 3 p\right)}
\biggl( e^{-c2^{2j} t} c_j \|R_0\|_{B^{-1+\frac 3 p}_{p,2}} + 
\int_{0}^t e^{-c2^{2j}(t-t')}
c_{j}(t')\|F(t')\|_{B^{-1+\frac 3 p}_{p,2}}dt'
\biggr)\\
& \leq  &  C  2^{j\left(1-\frac 3 p\right)}
\Biggl( e^{-c2^{2j} t} c_j \|R_0\|_{B^{-1+\frac 3 p}_{p,2}} + 
\biggl(\int_{0}^t e^{-c2^{2j}(t-t')}
\|F(t')\|_{B^{-1+\frac 3 p}_{p,2}}dt'\biggr)^{\frac 1 2}\\
&&\qquad\qquad\qquad\qquad\qquad\quad{}\times
\biggl(\int_{0}^t e^{-c2^{2j}(t-t')}
c^2_{j}(t')\|F(t')\|_{B^{-1+\frac 3 p}_{p,2}}dt'\biggr)^{\frac 1 2}\Biggr).
\eeno
Then we infer immediately that
\beno
\|\Delta_{j}\cR_{0,\lam}\|_{L^\infty(\R^+;L^p)}  & \lesssim &
 2^{j\left(1-\frac 3 p\right)}
 c_j \|R_0\|_{B^{-1+\frac 3 p}_{p,2}} \\
 & & \ {}+  2^{j\left(1-\frac 3 p\right)}
\biggl(\int_{0}^\infty
c^2_{j}(t)\|F(t)\|_{B^{-1+\frac 3 p}_{p,2}}dt\biggr)^{\frac 1 2}
\|F\|^{\frac 1 2}_{L^1(\R^+;B^{-1+\frac 3 p}_{p,2})}
\quad\hbox{and}\\
\|\Delta_{j}\cR_{0,\lam}\|_{L^2(\R^+;L^p)}   & \lesssim &  
2^{-j\frac 3 p}
 c_j \|R_0\|_{B^{-1+\frac 3 p}_{p,2}} \\
 & & \ {}+2^{-j\frac 3 p}
\biggl(\int_{0}^\infty
c^2_{j}(t)\|F(t)\|_{B^{-1+\frac 3 p}_{p,2}}dt\biggr)^{\frac 1 2}
\|F\|^{\frac 1 2}_{L^1(\R^+;B^{-1+\frac 3 p}_{p,2})}.
\eeno
This gives
\[
\|\cR_{0}\|_{X_{\lam}}\lesssim    \|R_0\|_{B^{-1+\frac 3 p}_{p,2}} +   
\|F\|_{L^1(\R^+;B^{-1+\frac 3 p}_{p,2})} .
\]
It follows that the smallness condition\refeq{smallnessPNS} implies 
 precisely condition~(\ref{estimF0}). So we can apply Lemma\refer{cacciopoli+} which
gives a global, unique solution~$R$ to~$(PNS)$ such that
\[
R\in L^\infty(\R^+;B^{-1+\frac 3 p}_{p,2})\cap
L^2 (\R^+;B^{\frac 3 p}_{p,2}).
\]
We leave the classical proof of the continuity in time to the
reader. Theorem~\ref{PNSglobal}  is proved, provided we prove Proposition\refer{propestimconjugaison}. \ef

\medskip

\noindent{\bf Proof of Proposition\refer{propestimconjugaison} }
It relies mainly on Lemma\refer{regulheatflow} and in a Bony type 
decomposition. In order 
to prove the estimate on~$B_{NS}$, let us observe that 
Lemma\refer{regulheatflow} implies that
\[
\|\D_{j}(B_{NS}(R,R'))_{\lam}(t)\|_{L^p} \leq
e^{\lam\int_{0}^{\infty}\|u^{(0)}(t)\|^2_{L^\infty}dt}
\int_{0}^t e^{-c2^{2j}(t-t')}
\|\D_{j}{\bf P}\dive (R_{\lam}(t')\otimes R'_{\lam}(t'))\|_{L^p}dt'.
\]
Proposition 3.1 of\ccite{chemin20} implies that
\[
\|\D_{j}{\bf P}\dive (R_{\lam}\otimes R'_{\lam})
\|_{L^2(\R^+;L^p)} \leq  C c_{j}2^{j\left(-1+\frac 3 p\right)}
\|R\|_{X_{\lam}}\|R'\|_{X_{\lam}}.
\]
Young's inequality in time  ensures the estimate on~$B_{NS}$. 

\medskip 

\noindent The study 
of~$L_{0}$ follows the ideas of\ccite{chemin20}.
 Let us decompose~$(L_{0}a)_{\lam}\eqdefa L_{0}a_{\lam}$ as a sum  of two 
operators~$L_{1,\lam}$ and~$L_{2,\lam}$ defined by
\beno
(L_{n,\lam}R)(t) & \eqdefa & \int_{0}^t 
e^{(t-t')\Delta-\lam\int_{t'}^t\|u^{(0)}(t'')\|^2_{L^\infty}dt''}
{\bf P}\dive\cT_{n}(u^{(0)}(t'),R_{\lam}(t'))dt'\quad\hbox{with}\\
\cT_{1}(a,b) & \eqdefa & \sum_{j} \left(\D_{j}a\otimes S_{j-1}b+
S_{j-1}b\otimes \D_{j}a\right)\quad\hbox{and}\\
\cT_{2}(a,b) & \eqdefa & \sum_{j} \left( S_{j+2}a\otimes \D_{j}b+
\D_{j}b\otimes S_{j+2}a\right).
\eeno
As
\[
\D_{j}\cT_{1}(a,b) =\sum_{|j'-j|\leq 5} \D_{j}
\left(\D_{j'}a\otimes S_{j'-1}b+S_{j'-1}b\otimes \D_{j'}a\right),
\]
we have
\[
\|\D_{j}\cT_{1}(a,b)\|_{L^p}\leq C \sum_{|j'-j|\leq 5}
\|\D_{j'}a\|_{L^\infty}\|S_{j'-1}b\|_{L^p}.
\]
Noticing that
\[
\|S_{j-1} R_\lambda\|_{L^\infty(\R^+;L^p)} \leq 
c_{j}2^{j\left(1-\frac 3 p\right)}
\|R\|_{X_{\lam}},
\]
we obtain
\[
\|\D_{j}\cT_{1}(u^{(0)}(t),R_{\lam}(t))\|_{L^p} \leq C
c_{j}2^{j\left(1-\frac 3 p\right)} \|R\|_{X_{\lam}}
\|u^{(0)}(t)\|_{L^\infty}.
\]
Using Bernstein's inequality and Lemma\refer{regulheatflow}, we have therefore
\beno
\|\D_{j}(L_{1,\lam}R)(t)\|_{L^p} & \leq & C 2^j\int_{0}^t 
e^{-c2^{2j}(t-t')-\lam\int_{t'}^t\|u^{(0)}(t'')\|^2_{L^\infty}dt''}
\|\cT_{1}(u^{(0)}(t'),R_{\lam}(t'))\|_{L^p}dt'\\
& \leq & C c_{j}2^{j\left(2-\frac 3 p\right)} \|R\|_{X_{\lam}}
\int_{0}^t e^{-c2^{2j}(t-t')-\lam\int_{t'}^t\|u^{(0)}(t'')\|^2_{L^\infty}dt''}
\|u^{(0)}(t')\|_{L^\infty}dt'.
\eeno
Thus we get, by Young's inequality, 
\beq
\label{demopropestimconjugaisoneq1}
\|\D_{j}(L_{1,\lam}R)\|_{L^\infty(\R^+;L^p)} 
+ 2^{j}\|\D_{j}(L_{1,\lam}R)\|_{L^2(\R^+;L^p)} 
\leq  \frac C {\lam^{\frac 1 2}}
 c_{j}2^{j\left(1-\frac 3 p\right)} \|R\|_{X_{\lam}}. 
\eeq
Let us now estimate~$ L_{2,\lam} R$. As
\[
\D_{j}\cT_{2}(a,b) =\sum_{j'-j\geq N_{0}} \D_{j}
\left(\D_{j'}a\otimes S_{j'-1}b+S_{j'-1}b\otimes \D_{j'}a\right),
\]
we have
\beno
\|\D_{j}\cT_{2}(a,b)\|_{L^p} & \leq  & C \sum_{j'-j\geq N_{0}} 
\|S_{j'+2}a\|_{L^\infty}\|\D_{j'}b\|_{L^p}\\
& \leq  & C \|a\|_{L^\infty}\sum_{j'-j\geq N_{0}} \|\D_{j'}b\|_{L^p}.
\eeno
As for the estimate of~$L_{1,\lam}$ we get that
\[
\longformule{\|\D_{j}(L_{2,\lam}R)(t)\|_{L^p} \leq   C2^j
\sum_{j'\geq j-N_{0}}\int_{0}^t e^{-c2^{2j}(t-t')-\lam\int_{t'}^t\|u^{(0)}(t'')\|^2_{L^\infty}dt''}
}
{
{}\times\|u^{(0)}(t')\|_{L^\infty}\|\D_{j'}R_{\lam}(t')\|_{L^p}dt'.
}
\]
The Cauchy-Schwarz inequality implies that
\[
\longformule{\|\D_{j}(L_{2,\lam}R)(t)\|_{L^p} \leq   C2^j
\sum_{j'\geq j-N_{0}}\biggl(\int_{0}^t e^{-c2^{2j}(t-t')}
\|\D_{j'}R_{\lam}(t')\|^2_{L^p}dt'\biggr)^{\frac 1 2}}
{
{}\times\biggl(\int_{0}^t 
e^{-2\lam\int_{t'}^t\|u^{(0)}(t'')\|^2_{L^\infty}dt''}
\|u^{(0)}(t')\|^2_{L^\infty}dt'\biggr)^{\frac 1 2}.
}
\]
Then we infer that
\[
\longformule{
2^{j\frac 3 p} \Bigl(\|\D_{j}(L_{2,\lam}R)\|_{L^\infty(\R^+;L^p)} 
+ 2^{j}\|\D_{j}(L_{2,\lam}R)\|_{L^2(\R^+;L^p)}\Bigr) 
}
{{}\leq 
\frac C  {\lam^{\frac 1 2}} \sum_{j'\geq j-N_{0}}
2^{(j-j')\frac 3 p} 2^{j'\frac 3 p}
\|\D_{j'}R_{\lam}\|_{L^2(\R^+;L^p)}.
}
\]
Young's inequality on series and\refeq{demopropestimconjugaisoneq1}
allow to conclude the proof of Proposition\refer{propestimconjugaison}.\ef

\section{End of the proof of Theorem~\ref{mainresult}}\label{proof}
\setcounter{equation}{0}
Now we are ready to prove Theorem~\ref{mainresult}. The idea, as
presented in the introduction, is to write
$$
u = u^{(0)} + R,
$$
where~$ R$ satisfies~$(PNS)$ with~$ R_0 = 0$ and~$ F =  - (\Id - \M)
{\bf P}
(u_{F}  \cdot \nabla u_{F}) 
- Q(u_{F} ,  u_{2D}) $, and where~$
u^{(0)} = u_{F} +  u_{2D}$.  
 According to the assumptions of Theorem~\ref{mainresult}, we know
 that~$ u_0$  belongs to~$\cI_p(A,B)$, so in particular
by~$ (\H3)$ we have
$$
\|F\|_{L^1(\R^+;B^{-1+\frac 3 p}_{p,2})} \leq B.
$$
Moreover by Corollary~\ref{u0L2linfty} we have
$$
\|u^{(0)}\|_{L^2 (\R^+;L^\infty)}^2 \lesssim A^2 (1 + A\log (e + A))^2 .
$$
Due to Theorem~\ref{PNSglobal}, the global wellposedness of~$ (PNS)$
is guaranteed if
$$
\|F\|_{L^1(\R^+;B^{-1+\frac3p}_{p,2})} \leq C_0^{-1} e^{-C_0
  \|u^{(0)}\|^2_{L^2(\R^+;L^\infty)}} .
$$
Clearly the smallness assumption~(\ref{smallness}) implies directly
that inequality, so  under the assumptions of
Theorem~\ref{mainresult}, we have
\[
R\in C_{b}(\R^{+};B^{-1+\frac3p}_{p,2}) 
\cap L^2(\R^+;B^{\frac3p}_{p,2}).
\] 
To end the proof of Theorem\refer{mainresult} we still need
to prove that~$ u$ is in~$ C_b(\R^+;H^\frac12) \cap
L^2(\R^+;H^\frac32)$. 
It is well known (see for instance\ccite{chemin10}) that the blow up 
condition for~$H^{\frac 1 2}(\TT^3)$ data is the blow up of the
norm~$L^2$ in time with values in~$H^{\frac  32}$. As~$u_{0}$ is 
in~$H^{\frac 1 2}$, so are~$\overline u_{0}$ and~$\wt u_{0}$.
Then thanks to the propagation of regularity in~$(NS2D)$
(see for instance\ccite{chemin10}) and the properties of the heat 
flow,~$u_{F}$ and~$u_{2D}$ belong to
 \[
L^\infty(\R^+;H^{\frac 1 2})\cap L^2(\R^+;H^{\frac  32}) 
\quad\hbox{and thus to}\quad
L^\infty(\R^+;B^{-1+\frac 3 p}_{p,2})\cap 
L^2(\R^+;B^{\frac 3 p}_{p,2})
\]
by the embedding recalled in~(\ref{embedding}). 
Thus as~$R$ belongs also to this space, it is enough to prove the following 
 blow up result, which we prove for  the reader's convenience. 
\begin{prop}
\label{trivialblowup}
{\sl If the maximal  time~$T^\star$ of existence in~$
L_{loc}^\infty(\R^+;H^{\frac 1 2})
\cap L^2_{loc}(\R^+;H^{\frac  32}) 
$
of a solution~$u$ of~$(NS)$  is finite, then for any~$p$,
\[
\int_{0}^{T^\star}\|u(t)\|^4_{B^{-\frac 1 2 +\frac 3 p}_{p,\infty}} dt
=+\infty.
\] }
\end{prop}

\noindent{\bf Proof of Proposition\refer{trivialblowup} }
An energy estimate in~$H^{\frac 1 2}$ gives, for some positive~$c$,
\[
\|u(t)\|_{H^{\frac 1 2}}^2
+c\int_{0}^{t} \|u(t')\|^2_{H^{\frac 3 2}}dt'\leq 
\|u_{0}\|_{H^{\frac 1 2}}^2 +
2 \int_{0}^{t} \left(\dive (u(t')\otimes u(t'))
|u(t')\right)_{H^{\frac 1 2}}dt'.
\]
We can assume that~$p>6$. Laws of product in Besov spaces imply  that
\[
\|u(t')\otimes u(t')\|_{H^{\frac 1 2}} \leq C 
\|u(t')\|_{B^{-\frac 1 2 +\frac 3 p}_{p,\infty}}
\|u(t')\|_{H^1}.
\]
Thus by interpolation we infer that
\[
\left(\dive (u(t')\otimes u(t'))
|u(t')\right)_{H^{\frac 1 2}} \leq C 
\|u(t')\|_{B^{-\frac 1 2 +\frac 3 p}_{p,\infty}}
\|u(t')\|^{\frac 1 2}_{H^{\frac 1 2}}
\|u(t')\|^{\frac 3 2}_{H^{\frac 3 2}}.
\]
Using the convexity inequality~$ab \leq 3/4 a^{\frac 4 3 }+ 1/4 b^4$
gives
\[
\|u(t)\|_{H^{\frac 1 2}}^2
+\frac c 2\int_{0}^{t} \|u(t')\|^2_{H^{\frac 3 2}}dt'\leq 
\|u_{0}\|_{H^{\frac 1 2}}^2 +
C \int_{0}^t \|u(t')\|^4_{B^{-\frac 1 2 +\frac 3 p}_{p,\infty}}
\|u(t')\|^{2}_{H^{\frac 1 2}} dt'.
\]
A Gronwall lemma concludes the proof of Proposition~\ref{trivialblowup}, and
therefore of Theorem\refer{mainresult}. \ef

\section{Proof of Theorem~\ref{example}} \label{exampleproof}
\setcounter{equation}{0}
In this final section we shall prove Theorem~\ref{example}. In order to do so, two points
must be checked: first, that the initial data defined in the statement of the theorem satisfies
the assumptions of Theorem~\ref{mainresult}, namely the nonlinear smallness
 assumption~(\ref{smallness}),
 in which case the global wellposedness will
follow as a consequence of that theorem. Second, that the initial data satisfies the lower
bound~(\ref{largeC-1}). Those two points are dealt with in Sections~\ref{checkassumption}
and~\ref{lowerbound} respectively.

\subsection{The nonlinear smallness assumption}
\label{checkassumption}
Let us check that the initial data defined in the statement of  Theorem~\ref{example} belongs to
 the space~$\cI_p(A,B)$ with the smallness
 condition~(\ref{smallness}).  Recall that~$A$ and~$B$ are chosen
so that
\beno
\displaystyle (\H1) \quad    &\ds \|\overline u_{0}\|_{L^2(\TT^2) } + 
\|\M  {\bf P} (  u_{F}\cdot \nabla   u_{F})\|_{L^1(\R^+;L^2(\TT^2))} \leq A & \\
\displaystyle (\H2)  \quad    &\ds \|\widetilde{u_{0}}\|_{
  B^{-1}_{\infty,2}} \leq A&\\
\displaystyle (\H3)  \quad        &\ds\left\|
(\Id - \M) {\bf P} (u_{F}\cdot \nabla  u_{F} )
  +Q  ( u_{2D} ,  u_{F} )
 \right\|_{L^1(\R^+; B^{-1+\frac3p}_{p,2})}   \leq B .&
\eeno
\medskip
\noindent Let us start with Assumption~$(\H1)$. We first  notice
 directly that~$\overline u_{0} = 0$,  so 
we just have to check that~$\M  {\bf P} (  u_{F}\cdot \nabla   u_{F})$
 belongs to~$L^1(\R^+;L^2(\TT^2))$, and to compute its bound.

\noindent We have
$$
 u_{F}\cdot \nabla   u_{F} = \dive (u_{F} \otimes u_{F})
\quad\hbox{hence}\quad
\M   ( u_{F}\cdot \nabla   u^j_{F})  =  \M \dive_{h}  (u^j_{F}
 u^h_{F})  \quad \mbox{for} \quad j \in \{1,2,3\}.
$$
On the one hand, we have
\begin{eqnarray*}
\dive_{h}  (u^3_{F}u^h_{F})(x)  & = & 
 N \dive_{h} \left(-\dive_{h}(e^{t \Delta_{h}} v_{0}^h)
 e^{t \Delta_{h}} v_{0}^h\right)(x_{h})
\left( e^{t \partial_{3}^{2}}\sin(Nx_{3})\right)
\left(e^{t \partial_{3}^{2}}\cos(Nx_{3})\right) \\
& = & \frac N2 e^{-2tN^{2}} 
\dive_{h} \left(-\dive_{h}(e^{t \Delta_{h}} v_{0}^h)
 e^{t \Delta_{h}} v_{0}^h\right)(x_{h})\sin(2Nx_{3}),
\end{eqnarray*}
which  implies that~$\M   ( u_{F}\cdot \nabla   u^3_{F})  = 0$.
Notice that in particular, since~$ \overline u_0 = 0$, we infer that
\beq\label{u2d30}
\forall t \geq 0, \quad u_{2D}^3 (t) = 0.
\eeq
On the other hand, we have
\beno
 \dive_{h}(u^h_{F} \otimes u_{F}^{h}) (x) & = & 
 N^2 \dive_{h}\left(e^{t \Delta_{h}} v_{0}^h\otimes
 e^{t \Delta_{h}} v_{0}^h\right)(x_{h})
 \left(e^{t \partial_{3}^{2}}\cos(Nx_{3})\right)^2\\
 & = & \frac {N^{2}}2  e^{-2tN^{2}}
 \dive_{h}\left(e^{t \Delta_{h}} v_{0}^h\otimes
 e^{t \Delta_{h}} v_{0}^h\right)(x_{h}) (1+\cos (2Nx_{3})).
\eeno
Using the frequency localization of~$ v_{0}^h$ and Bernstein's
 inequality~(\ref{bernstein}), we get
\begin{eqnarray*}
\left\| \M (u_{F}\cdot \nabla u_{F}^{h})  \right\|_{L^{2}(\TT^{2})} & \leq&
\frac {N^{2}}2  e^{-2tN^{2}}N_{0} \|e^{t \Delta_{h}} v_{0}^h\|_{L^{4}(\TT^{2})}^{2 } \\
 & \leq& C N_{0}^2 N^{2} e^{-2tN^{2}} \| v_{0}^h\|_{L^{2}(\TT^{2})}^{2 }.
\end{eqnarray*}
 Finally we infer that
\ben
\left\|\M (u_{F} \cdot \nabla u_{F}^{h}) 
\right\|_{L^{1}(\R^{+};L^{2}(\TT^{2}))} 
& \leq  & C N_{0}^2  \| v_{0}^h\|_{L^{2}(\TT^{2})}^{2 }\\
\label{eqH1}& \leq  & C_{N_{0}} (\log N)^{\frac 2 9}\nonumber.
\een

\bigskip

\noindent Let us now consider  Assumption~$(\H2)$. Since~$\overline u_{0} = 0$, it simply consists
in computing the~$ B^{-1}_{\infty,2}$ norm of~$u_{0}$. We have
$$
u_{0}^{h} (x)= N v_{0}^h(x_{h}) \cos (Nx_{3}) ,
$$
and by definition of Besov norms,
$$
\|u_{0}^{h}\|_{B^{-1}_{\infty,2}} = \left\|
\tau^{\frac12} \|e^{\tau \Delta} u_{0}^{h}\|_{L^{\infty}}
 \right\|_{L^{2}(\R^{+},\frac{d\tau}{\tau})}=
  \|e^{\tau \Delta} u_{0}^{h}\|_{L^2(\R^+;L^{\infty})} .
$$
It is easy to see that
\begin{eqnarray*}
\|e^{\tau \Delta} u_{0}^{h}\|_{L^{\infty}} & = &  
N \| e^{\tau \Delta_{h}}v_{0}^h(x_{h})
 e^{\tau \partial_{3}^{2}} \cos (Nx_{3}) \|_{L^{\infty}} \\
&\leq & C {N_{0}} N e^{-\tau N^{2}} \|v_{0}^h\|_{L^{2}}.
\end{eqnarray*}
It follows that
\begin{eqnarray*}
\|u_{0}^{h}\|_{B^{-1}_{\infty,2}} &\leq & C {N_{0}} 
N  \|v_{0}^h\|_{L^{2}}  
\| e^{-\tau N^{2}} \|_{L^{2}(\R^{+})}\\
 &\leq & C {N_{0}}\|v_{0}^h\|_{L^{2}}  .
\end{eqnarray*}

\noindent The computation is similar for~$u_{0}^{3}$,  so we get, for~$N$ 
large enough,
\beq
\label{eqH2}
\|u_{0}\|_{B^{-1}_{\infty,2}} \leq C_{N_{0}}(\log N)^{\frac 1 9}.
\eeq 
Thus one   can choose
for the parameter~$A$ in~$(\H1)$ and~$(\H2)$
\begin{equation}\label{defA}
A =  C _{N_{0}}  (\log N)^{\frac 2 9}.
\end{equation}

\bigskip 

\noindent Finally let us consider Assumption~$(\H3)$. We shall start
with~$ (\Id - \M){\bf P} (u_F \cdot \nabla u_F)$. We have
$$
 (\Id - \M)  (u_F \cdot \nabla u_F) =  (\Id - \M)  (u_F^h \cdot \nabla^h
 u_F) +  (\Id - \M)  (u_F^3 \partial_3 u_F) ,
$$
and we shall concentrate on the first term, as both are treated in the
same way. We compute
\beno
 (\Id - \M)  (u_F^h \cdot \nabla^h u_F^h) (x) & = & \frac{N^2}{2}   
 \Bigl(e^{t\Delta_{h}} v_{0}^h \cdot
 \nabla^h   e^{t\Delta_{h}} v_{0}^h \Bigr)(x_{h}) e^{-2tN^2} \cos (2Nx_3)
\quad \mbox{and} \\
 (\Id - \M)  (u_F^h \cdot \nabla^h u_F^3) & = &
  -\frac{N}{2} \Bigl ( e^{t\Delta_{h}} v_{0}^h \cdot
 \nabla^h   e^{t\Delta_{h}} \dive_h v_{0}^h \Bigr)(x_{h})  
 e^{-2tN^2} \sin  (2Nx_3) .
\eeno
So
\beno
\| (\Id - \M) (u_F^h \cdot \nabla^h
 u_F^h) \|_{B^{-1+\frac3p}_{p,2}} & = & 
\left\|
\tau^{\frac12 - \frac3{2p} } \|e^{\tau\Delta}  (\Id - \M) (u_F^h \cdot \nabla^h
 u_F^h)\|_{L^p}
\right\|_{L^2(\R^+,\frac{d\tau}{\tau})}
\\
&\leq & C_{N_0}\frac{N^2}{2}  e^{-2tN^2} \left\|
\tau^{\frac12 - \frac3{2p} }  e^{-2\tau N^2}
\right\|_{L^2(\R^+,\frac{d\tau}{\tau})} \| v_{0}^h\|_{L^2}^2 \\
&\leq & C_{N_0}N^2 e^{-2tN^2} N^{\frac3p -1}  \| v_{0}^h\|_{L^2}^2 .
\eeno
It follows that
\beq\label{one}
\| (\Id - \M) (u_F^h \cdot \nabla^h
 u_F^h) \|_{L^1(\R^+;B^{-1+\frac3p}_{p,2})} \leq  C_{N_0} N^{\frac3p -1}  \| v_{0}^h\|_{L^2}^2 ,
\eeq
and similarly
\beq\label{two}
\| (\Id - \M) (u_F^h \cdot \nabla^h
 u_F^3) \|_{L^1(\R^+;B^{-1+\frac3p}_{p,2})} \leq  C_{N_0} N^{\frac3p -2}  \| v_{0}^h\|_{L^2}^2 .
\eeq
Finally let us estimate the term~$ Q(u_{2D},u_F)$. Since
by~(\ref{u2d30}), $u_{2D}^3 $ is identically equal to zero, we have
$$
Q(u_{2D},u_F) = {\bf P} \dive_h (u_{2D}  \otimes u_F +  u_F \otimes u_{2D} )
$$
so
$$
\|Q(u_{2D},u_F^h)\|_{L^1(\R^+;B^{-1+\frac3p}_{p,2})} \leq N
\|e^{-tN^2}\dive_h (e^{t\Delta_h}
v_{0}^h \otimes u_{2D}^h)(x_{h}) \cos(Nx_3)\|_{L^1(\R^+;B^{-1+\frac3p}_{p,2})} .
$$
we shall only compute that term, as~$Q(u_{2D},u_F^3) $ is estimated
similarly (and contributes in fact one power less in~$ N$). 
Sobolev embeddings imply that~$ H^s (\TT^2)\hookrightarrow L^p(\TT^2)$
for~$ s  \eqdefa 1-\frac2p\cdotp$~So
\beno
\|\dive_h (e^{t\Delta_h}
v_{0}^h \otimes u_{2D}^h) \|_{L^p} &\leq& \|e^{t\Delta_h}
v_{0}^h \cdot \nabla^h  u_{2D}^h\|_{L^p} + \| u_{2D}^h \cdot \nabla^h e^{t\Delta_h}
v_{0}^h \|_{L^p} \\
&\leq& \|e^{t\Delta_h} v_{0}^h\|_{L^\infty} \| \nabla^h
u_{2D}^h\|_{L^p} + \| u_{2D}^h\|_{L^p}  \|e^{t\Delta_h} v_{0}^h\|_{L^\infty} \\
&\leq& C N_0  \|v_{0}^h\|_{L^2} \| u_{2D}^h\|_{H^{s+1}} + C N_0^2
\|v_{0}^h\|_{L^2}  \| u_{2D}^h\|_{H^{s}}.
\eeno
Propagation of regularity for the two dimensional Navier-Stokes
equations is expressed by
\[
\|u_{2D}\|_{L^\infty(\R^+;H^2)} \leq 
\|\M(u_{F}\cdot\nabla u_{F})\|_{L^1(\R^+;H^2)} 
e^{C\|u_{2D}\|_{L^2(\R^+;H^1)}^2}.
\]
Using\refeq{eqH1} and  that
the Fourier transform of~$\M(u_{F}\cdot\nabla u_{F})$ is supported
in~$[-2N_{0},2N_{0}]^2$, we~get
\beno
\|u_{2D}\|_{L^\infty(\R^+;H^2)}  & \leq  & CN_{0}^2
\|\M(u_{F}\cdot\nabla u_{F})\|_{L^1(\R^+;L^2)} 
e^{C_{N_{0}}\|v_{0}^h\|^4_{L^2}} \\
 & \leq  &  C_{N_{0}}\|v_{0}^h\|^2_{L^2}
e^{C_{N_{0}}\|v_{0}^h\|^4_{L^2}}.
\eeno
Therefore we obtain
\beno
\| Q(u_{2D},u_F)(t)\|_{B^{-1+\frac3p}_{p,2}} &\leq &  
C_{N_0}  N e^{-tN^2}
\left\|\tau^{\frac12 - \frac3{2p}}e^{-\tau N^2}
\right\|_{L^2(\R^+;\frac{d\tau}{\tau})} 
\|v_{0}^h\|^3_{L^2} e^{C_{N_{0}}\|v_{0}^h\|^4_{L^2}} \\
 &\leq &  C_{N_0} N^{\frac3 p} e^{-tN^2} \|v_{0}^h\|^3_{L^2} 
 e^{C_{N_{0}}\|v_{0}^h\|^4_{L^2}}.
\eeno
Finally
\[
\| Q(u_{2D},u_F)\|_{L^1(\R^+;B^{-1+\frac3p}_{p,2})} \leq  C_{N_0}
 N^{\frac3 p -2} \|v_{0}^h\|^3_{L^2} 
 e^{C_{N_{0}}\|v_{0}^h\|^4_{L^2}}.
 \]
Together with\refeq{one} and\refeq{two}, this gives
\[
\longformule{
\left\|(\Id - \M) {\bf P} (u_{F}\cdot \nabla  u_{F} )
  +Q  ( u_{2D} ,  u_{F} )
 \right\|_{L^1(\R^+; B^{-1+\frac3p}_{p,2})} 
 }  
{ {}\leq
C_{N_0} N^{\frac3 p -1}  \|v_{0}^h\|^2_{L^2} 
\Bigl(1 +N^{-1}\|v_{0}^h\|_{L^2} e^{C_{N_{0}}\|v_{0}^h\|^4_{L^2}}
\Bigr) . 
}
\]
Using that~$\|v_{0}^h\|_{L^2(\TT^2)}\leq (\log N)^\frac 1 9$, we infer that,
for~$N$ large enough,
\[
\left\|(\Id - \M) {\bf P} (u_{F}\cdot \nabla  u_{F} )
  +Q  ( u_{2D} ,  u_{F} ) \right\|_{L^1(\R^+; B^{-1+\frac3p}_{p,2})} 
\leq C_{N_{0}} N^{\frac 3 p-1} (\log N)^{\frac 2 9}.
\]
Choosing~$p\geq 6$ gives, still for~$N$ large enough,
\[
\left\|(\Id - \M) {\bf P} (u_{F}\cdot \nabla  u_{F} )
  +Q  ( u_{2D} ,  u_{F} ) \right\|_{L^1(\R^+; B^{-1+\frac3p}_{p,2})} 
\leq  N^{-\frac 1 4}.
\]
We can therefore choose for the parameter~$ B$ in~$ (\H3)$ 
the value~$B=N^{-\frac 1 4}$.
Let us check that with such choices of~$A$ and~$ B$, the smallness
assumption\refeq{smallness} holds. 
With the choice of~$A=C_{N_{0}}(\log N)^{\frac 2 9}$
made in\refeq{defA}, we have, for~$N$ large enough,
\beno
\exp \left(C_{0}A^2(1+A\log A)^{ 2}\right) & \leq & \exp \left(C_{N_{0}}
(\log^{\frac 8 9}N)(\log \log N)\right)\\
 & \leq  & \exp \Bigl(\frac 1 8 \log N\Bigr)\\
 & \leq & N^{\frac 1 8}.
\eeno
Since~$B=N^{-\frac 1 4}$, the smallness assumption\refeq{smallness} is guaranteed
for large enough~$N$, and Theorem~\ref{mainresult} yields the global wellposedness of the
system with that initial data. 

\subsection{The lower bound} \label{lowerbound}
Let us now check that the initial data~$u_{0}^{h}$ satisfies the lower
bound\refeq{largeC-1}. 
We recall that the~$ B^{-1}_{\infty,\infty}$ norm is defined by
$$
\|u_{0}^{h}\|_{ B^{-1}_{\infty,\infty}} = \sup_{t \geq 0} t^{\frac12}
\|e^{t\Delta} u_{0}^{h}\|_{L^{\infty}(\TT^{3})}.
$$
An easy computation, using the explicit formulation of~$u_{0}^{h}$, enables us to write
that
\begin{eqnarray*}
e^{t\Delta} u_{0}^{h} (x)& = & N e^{t\Delta_{h}} v_{0}^h (x_{h})
e^{t\partial_{3}^{2}} \cos (Nx_{3})\\
& = & N  e^{t\Delta_{h}} v_{0}^h (x_{h}) e^{-t N^{2}}  \cos (Nx_{3}).
\end{eqnarray*}
It follows that
\begin{eqnarray*}
\|e^{t\Delta} u_{0}^{h}\|_{L^{\infty}(\TT^{3})} & =&  N e^{-t N^{2}} 
 \| e^{t\Delta_{h}} v_{0}^h \|_{L^{\infty}(\TT^{2})}\\
& \geq& \frac N{2\pi} e^{-t N^{2}}\| e^{t\Delta_{h}} 
v_{0}^h \|_{L^{2}(\TT^{2})}\\
& \geq& \frac N{2\pi} e^{-2t N^{2}} \|v_{0}^h \|_{L^{2}(\TT^{2})},
\end{eqnarray*}
for~$N \geq N_{0}$, using the fact that the frequencies of~$v_{0}^h $ are smaller than~$N_{0}$. 
Finally we have
\begin{eqnarray*}
\|u_{0}^{h}\|_{ B^{-1}_{\infty,\infty}}  & \geq& 
  \frac N{2\pi}\|v_{0}^h \|_{L^{2}(\TT^{2})}
\sup_{t \geq 0}\left (t^{\frac12} e^{-2t N^{2}}\right) \\
  & \geq& \frac{1}{4 \pi \sqrt e} \|v_{0}^h \|_{L^{2}(\TT^{2})},
\end{eqnarray*}
and Theorem\refer{example} follows. \ef

\end{document}